\documentclass[a4]{amsart}
\usepackage{amsmath,amssymb,amsthm,enumerate,epsfig,subfigure}
\usepackage[monochrome]{color}

\newcommand{\Q}{{\mathbb Q}}
\newcommand{\Z}{{\mathbb Z}}

\newcommand{\N}{{\mathbb N}}

\newtheorem{thm}{Theorem}

\newtheorem{cor}{Corollary}
\newtheorem{conj}{Conjecture}

\newcommand{\Mex}{\mathrm{Mex}}
\newcommand{\Gap}{\mathrm{Gap}}
\newcommand{\Bl}{\mathrm{Block}}

\begin{document}
\title{3 dimensional Wythoff Nim}
\author{Shigeki Akiyama}
\maketitle

\begin{abstract}
We introduce a new generalization of Wythoff Nim using
 three piles of stones. We show
that its P-positions have
finite difference properties and produce a partition of positive integers. 
Further, we give a 
conjecture that the P-positions approximate a half-line whose 
slope is described by algebraic numbers of degree 5.
\end{abstract}

\section{Introduction}
\label{Intro}
``Tsyanshidzi"\footnote{In the Japanese book 
`World of Mathematical Games' by S.~Hitotsumatsu, we find a claim that
this name came from a Russian translation of an ancient Chinese 
``Select stones''.
%'IÎFBo Tan told me that it may be EÎŽq Jianshizi
}
is an ancient Chinese game between two players. Its
winning strategy is first
studied by Wythoff \cite{Wythoff}, and it is called Wythoff Nim.
The rule is

\begin{itemize}
\item Given two piles of stones.
Each player may remove an arbitrary
number of stones from one of the piles or may remove an equal number of stones from both piles. The one who emptied both piles is a winner.
\end{itemize}

Such a complete information game that finishes in finite steps 
without a draw, there are 
relatively rare positions where the previous player has a winning 
strategy, which are called P-positions.

Let $\N_0$ be the set of non-negative integers and $\N=\N_0\setminus \{0\}$.
Denote the set of P-positions 
$\Pi:=\{(a_n,b_n)\in \N^2\ 
|\ a_n\le b_n,\ a_{n}<a_{n+1}, \ n=1,2,\dots \}$ of Wythoff Nim and 
define $\Mex(F)=\min(\N_0\setminus F)$. The set
$\Pi$ has several remarkable properties:

\begin{enumerate}
\item {\bf Mex Algorithm:} There exists a fast algorithm to compute 
$\Pi$ using $\Mex$ function.
\item {\bf Finite difference:}
 $\{ (a_{n+1}-a_n,b_{n+1}-b_n)\ |\ n\in \N\}$ is finite. Indeed it is 
$\{(1,2),(2,3)\}$.
\item {\bf Partition of $\N$:} $\{a_n\ |\ n\in \N\} \cup \{b_n\ |\ n\in \N\}$ forms a partition of $\N$.
\item {\bf Substitutive structure:} The difference sequence 
$(a_{n+1}-a_n)_{n=1,2,\dots}$ is the fixed point of the Fibonacci substitution $2\rightarrow 21,\ 1\rightarrow 2$.
\item {\bf Linear structure:}
$\Pi$ approximates a half-line emanating from the origin. More precisely
$$
a_n=\lfloor \omega n \rfloor,\ b_n=\lfloor \omega^2 n \rfloor
$$
with $\omega=(1+\sqrt{5})/2$.
\end{enumerate}

Because $\omega$ is irrational, Wythoff Nim looks 
very different and difficult from other well-studied 
combinatorial Nims, e.g., we do not know a quick
way to compute Grundy numbers. 
This game attracted many mathematicians because of the beauty of the P-positions and 
there are a lot of variants, trying to extend those properties, see \cite{WythoffV, Rigo_Game}.

In this paper, we introduce Wythoff 3D $L^1$-Nim:
A subtraction game on three piles of stones, which shares many properties 
of the original Wythoff Nim. 
A move $$\N_0^3\ni (a,b,c)\to (a-x,b-y,c-z)\in \N_0^3$$ is defined by
$$
\{(x,y,z)\in \Z^3\ |\ xyz(x-y)(y-z)(z-x)=0,\ x+y+z>0\}.
$$
Since $$|a-x|+|b-y|+|c-z|<|a|+|b|+|c|,$$ 
each move must decrease its $L^1$ norm, and we have the finite ending property.
Note that some of $x,y,z$ can be negative, i.e., one can increase the number of stones.
\bigskip

Here we list first several P-positions $(a,b,c)$ with $a\le b\le c$.
\bigskip

(0, 0, 0), (1, 2, 3), (4, 7, 10), (5, 9, 13), (6, 11, 16), (8, 15, 22), 
 (12, 21, 30), (14, 25, 36), (17, 29, 41), (18, 31, 44), (19, 34, 49), 
 (20, 37, 53), (23, 42, 61), (24, 45, 65), (26, 51, 74), (27, 54, 81), 
 (28, 56, 84), (32, 63, 92), (33, 67, 99), (35, 70, 105), (38, 75, 111), 
 (39, 78, 117), (40, 80, 120), (43, 85, 127), (46, 89, 132), (47, 91, 135), 
 (48, 93, 138), (50, 96, 142), \textcolor{blue}{(52, 101, 148)}, (55, 106, 156), 
 (57, 109, 161), \textcolor{blue}{(58, 113, 166)}, (59, 116, 173), (60, 118, 176), 
 (62, 121, 180), (64, 125, 186), (66, 128, 190), (68, 131, 194), 
 (69, 133, 197), (71, 136, 201).
\bigskip

whose differences $(b-a,c-b)$ are:
\bigskip

(0, 0), (1, 1), (3, 3), (4, 4), (5, 5), (7, 7), (9, 9), (11, 
  11), (12, 12), (13, 13), (15, 15), (17, 16), (19, 19), (21, 
  20), (25, 23), (27, 27), (28, 28), (31, 29), (34, 32), (35, 
  35), (37, 36), (39, 39), (40, 40), (42, 42), (43, 43), (44, 
  44), (45, 45), (46, 46),\textcolor{blue}{ (49, 47)}, (51, 50), (52, 52), \textcolor{blue}{(55, 
  53)}, (57, 57), (58, 58), (59, 59), (61, 61), (62, 62), (63, 
  63), (64, 64), (65, 65).
\bigskip

Readers may feel more natural to 
define the set of possible moves $(a,b,c)\to (a-x,b-y,c-z)$ by
$$
\{(x,y,z)\in \N_0^3\ |\ xyz(x-y)(y-z)(z-x)=0,\ x+y+z>0\},
$$
which we call 3D Wythoff positive Nim, i.e., we can only subtract non-negative
number of stones.

P-positions for two Nims are the same until (20, 37, 53) but the next one
is (23, 39, 58) for positive Nim while (23, 42, 61) for $L^1$-Nim.
The differences
$$(37-20, 53-37)=(17, 16),\quad (39-23, 58-39)=(16, 19)$$ 
share the number $16$. However the move
$$
(23, 39, 58) \rightarrow (20, 37, 53)
$$
is not allowed by positive Nim, and it is possible by the $L^1$ rule.
Therefore (23, 39, 58) is not a P-position in $L^1$-Nim.
In this manner for positive Nim, it is not possible to decide the next P-position only
by the occurrences of numbers 
and differences in previous P-positions. 
This makes this Nim difficult to analyze. 
We have not proved nor denied
any of the above five properties in the positive Nim.

In contrast, $L^1$-Nim is much simpler. 
Since some of $x,y,z$ could be negative,
the game rule sounds artificial but theoretically, it is very nice. 

In section \ref{Mex}, Theorem \ref{Algo} gives the property (1):
 a fast algorithm using $\Mex$ function. 
We show the property (3): 
the coordinates of P-positions give a partition of $\N$ 
in Corollary \ref{Par}.

Section \ref{Dynamics} is the technical key of this paper. 
We find that transitions of the game state $n\to n+1$
is understood as a discrete dynamics acting on three sets $X_n,
Y_n, Z_n$ of integers (see the end of section \ref{Dynamics} for their
definition).
We shall show that
$X_n$ receives $v_2(n+1)$, $Y_n$ receives $v_3(n+1)$ and
$Z_n$ receives $v_3(n+1)-v_1(n+1)$. 
$Y_n$ and $Z_n$ keeps their gaps
in this process.

As a consequence in section \ref{FD}, we show that
$X_n$ does not contain
four consecutive integers. This implies 
Theorem \ref{Run} and Corollary \ref{FinDiff}, which 
show that the
set of differences of adjacent P-positions are finite, 
i.e., the property (2). Indeed, the set has 56 differences and a simple structure, see Figure \ref{Prism}.
Further Corollary \ref{Balance} shows that each  P-position is very
close to the arithmetic progression of length three. 

Probably we cannot expect the property (4), 
we do not find a substitutive 
structure in the set of P-positions. 
However 
the above dynamics 
allow us to give a natural
conjecture of the approximate
linear structure of P-positions, that is, the property (5) in 
section \ref{LS}. 
The slope of the half-line is expected to be described by 
algebraic numbers of degree $5$, see Theorem \ref{Line}.

Klein-Fraenkel \cite{KleinFraenkel} introduced another higher dimensional
generalization of Wythoff Nim keeping the property (3).
In dimension three, it is called Wytlex(3), 
whose move is defined by
$$
\{(x,y,z)\in \Z^3\ |\ xyz(x-y)(y-z)(z-x)=0 \land
  (a-x,b-y,c-z)\ll (a,b,c)
\}
$$
where $\ll$ is the lexicographical order.
This game is very different from our $L^1$-Nim, e.g.,  
$$
(0,10,30) \rightarrow (0,9,40)
$$
is allowed in Wytlex(3) but its converse
$$
(0,9,40) \rightarrow (0,10,30)
$$
is allowed
in $L^1$-Nim. However it is interesting to note that
Theorem \ref{Algo} below proves that the sets of P-positions are exactly 
the same, which is numerically
conjectured in \cite{KleinFraenkel}.
Indeed the algorithms to produce $P$-positions are equivalent, 
see \cite[Theorem 1 and Lemma 3]{KleinFraenkel}.
Therefore all results in this paper apply to Wytlex(3) as well.

\section{Mex Algorithm and its consequences}
\label{Mex}

We start with the fast algorithm:

\begin{itemize}
\item $G_0=\{(0,0,0)\}$ and $n\leftarrow 0$
\item[$\dag$]$F_n=\{ x_i\ |\ x=(x_1,x_2,x_3)\in G_n,\  i\in \{1,2,3\} \}$
\item $v_1(n+1)=\Mex(F_n)$
\item $D_n=\bigcup_{x\in G_n} \{x_2-x_1,x_3-x_2,x_3-x_1\}$
\item $v_2(n+1)=\Mex \left((v_1(n+1)+D_n) \cup \{1,2,\dots,v_1(n+1)\} \cup F_n\right)$
\item $v_3(n+1)=\Mex \left((v_2(n+1)+D_n) \cup \{1,2,\dots,v_2(n+1)\} \cup F_n\right)$
\item $G_{n+1}=G_n \cup \{(v_1(n+1),v_2(n+1),v_3(n+1))\}$
\item $n\leftarrow n+1$ and go back to $\dag$.
\end{itemize}

\begin{thm}
\label{Algo}
$G=\bigcup_{n=0}^{\infty} G_n$ is the set of P-positions of $L^1$-Nim.
\end{thm}

\begin{proof}
Put $v(n)=(v_1(n),v_2(n),v_3(n))$. We prove by induction. 
Assume $G_{n}$ is the subset of P-positions of $L^1$-Nim.
Since $\{v_1(n+1),v_2(n+1), v_3(n+1)\}\cap F_n=\emptyset$,
if there is a move from $v(n+1)$ to the element of $G_n$, then
we must have either $v_2(n+1)-v_1(n+1)\in D_n$
or $v_3(n+1)-v_2(n+1)\in D_n$ or $v_3(n+1)-v_1(n+1)\in D_n$.
However, by definition of $v(n)$, the first
two inclusions are impossible.
The last inclusion $v_3(n+1)-v_1(n+1)\in D_n$ is also impossible because
$$
v_3(n+1)-v_1(n+1)>v_3(n)-v_1(n).
$$
The proof of this inequality comes later\footnote{
The proof of Theorem \ref{GAP} uses only the algorithm of $v(n)$.}
in Theorem \ref{GAP}.
Thus $v(n+1)$ is a possible $P$-position.

If $v(n+1)$ is not the next
$L^1$ minimum $P$-position
of this Nim. 
Then there exists another P-position
$(h_1,h_2,h_3)$ with $h_1+h_2+h_3\le v_{1}(n+1)+v_2(n+1)+v_3(n+1)$
and $h_1+h_2+h_3> v_{1}(n)+v_2(n)+v_3(n)$ 
by induction hypothesis. Since $v(n+1)$ 
is lexicographically the smallest possible
P-position, $h_1\ge v_1(n+1)$ must hold.  
Thus we have either $h_2<v_2(n+1)$ or $h_3<v_3(n+1)$. 

{\bf Case 1 $h_2<v_2(n+1)$}.
We claim that $h_2\in F_{n}$ or $h_2-h_1\in D_{n}$ holds.
Indeed, $v_2(n+1)$ is the minimum integer not less than
$v_1(n+1)+\Mex(D_{n})$ 
which is not in $F_{n}$. If $h_2-h_1<\Mex(D_{n})$ then
$h_2-h_1\in D_{n}$. If not, then we have 
$$v_1(n+1)+\Mex(D_{n})\le h_1+\Mex(D_{n})\le h_2
<v_2(n+1).$$
\noindent
This implies $h_2\in F_{n}$.
The claim is proved.
If $h_2\in F_{n}$ or $h_2-h_1\in D_{n}$ holds, then since 
$h_1+h_2+h_3>v_1(k)+v_2(k)+v_3(k)$ for $k\le n$, there exists a move
$(h_1,h_2,h_3)\to (v_1(k),v_2(k),v_3(k))$.
Thus $(h_1,h_2,h_3)$ is not a P-position.

{\bf Case 2 $h_2\ge v_2(n+1)$ and $h_3<v_3(n+1)$}.
We claim that $h_3\in F_{n}$ or $h_3-h_2\in D_{n}$
holds.
Indeed, $v_3(n+1)$ is the minimum integer not less than
$v_2(n+1)+\Mex(D_{n})$ which is not in $F_n$.
If $h_3-h_2< \Mex(D_n)$, then $h_3-h_2\in D_n$. 
If not, then
$$v_2(n+1)+\Mex(D_{n})\le h_2+\Mex(D_{n})\le h_3 
<v_3(n+1).$$
\noindent
This implies $h_3\in F_{n}$.
If  $h_3\in F_{n}$ or $h_3-h_2\in D_{n}$ holds,  then since 
$h_1+h_2+h_3>v_1(k)+v_2(k)+v_3(k)$ for $k\le n$, there exists a move
$(h_1,h_2,h_3)\to (v_1(k),v_2(k),v_3(k))$.
Thus $(h_1,h_2,h_3)$ is not a P-position. 

We have proved that 
$v(n+1)$ is the next
$L^1$ minimum $P$-position of this Nim, which is 
the next lexicographical $P$-position as well.
\end{proof}

\begin{cor}
\label{Par}
$$
 \{v_1(n)\ |\ n\in \N\} \sqcup
 \{v_2(n)\ |\ n\in \N\} \sqcup
 \{v_3(n)\ |\ n\in \N\} 
$$
is a partition of $\N$. 
\end{cor}

\begin{proof}
By the definition of $\Mex$, the union must be $\N$. 
Since $v_i(m)\in F_m$, $n>m$ implies 
$v_j(n)\neq v_i(m)$.
\end{proof}

\section{Dynamics on the set of integers}
\label{Dynamics}

The following example enlightens the idea of the inductive
proofs of Theorem \ref{GAP} and \ref{Fin}.
\begin{align*}
(\Mex F_{11}, \Mex D_{11})&=(23,19)\\
F_{11} \cap [\Mex F_{11},\Mex F_{11} + \Mex D_{11})&=\{25, 29, 30, 31, 34, 36, 37, 41\}\\
F_{11} \cap [\Mex F_{11} + \Mex D_{11},\infty)&=\{44, 49, 53\}\\
D_{11}\cap [\Mex D_{11}, \infty)&=\{22, 24, 26, 30, 33\}.
\end{align*}
We obtain $(v_1(12),v_2(12),v_3(12))=(23, {\bf 42, 61})$ with
differences $\{19, {\bf 38}\}$, and 
\begin{align*}
(\Mex F_{12}, \Mex D_{12})&=(24,20)\\
F_{12} \cap [\Mex F_{12},\Mex F_{12} + \Mex D_{12})&=\{25, 29, 30, 31, 34, 36, 37, 41, {\bf 42}\}\\
F_{12} \cap [\Mex F_{12} + \Mex D_{12},\infty)&=\{44, 49, 53, {\bf 61}\}\\
D_{12}\cap [\Mex D_{12}, \infty)&=\{22, 24, 26, 30, 33, {\bf 38}\}.
\end{align*}
We obtain $(v_1(13),v_2(13),v_3(13))=(24, {\bf 45, 65})$ with
differences $\{20, 21, {\bf 41}\}$, and 
\begin{align*}
(\Mex F_{13}, \Mex D_{13})&=(26,23)\\
F_{13} \cap [\Mex F_{13},\Mex F_{13} + \Mex D_{13})&=
\{29, 30, 31, 34, 36, 37, 41, 42, 44, {\bf 45}\}\\
F_{13} \cap [\Mex F_{13} + \Mex D_{13},\infty)&=
\{49, 53, 61, {\bf 65}\}\\
D_{13}\cap [\Mex D_{13}, \infty)&=
\{24, 26, 30, 33, 38, {\bf 41}\}.
\end{align*}
We have $(v_1(14),v_2(14),v_3(14))=(26, {\bf 51, 74})$ with
differences $\{23, 25, {\bf 48}\}$, and 
\begin{align*}
(\Mex F_{14}, \Mex D_{14})&=(27,27)\\
F_{14} \cap [\Mex F_{14},\Mex F_{14} + \Mex D_{14})&=
\{29, 30, 31, 34, 36, 37, 41, 42, 44, 45, 49, {\bf 51}, 53\}\\
F_{14} \cap [\Mex F_{14} + \Mex D_{14},\infty)&=
\{61, 65, {\bf 74}\}\\
D_{14}\cap [\Mex D_{14}, \infty)&=
\{30, 33, 38, 41, {\bf 48}\}.
\end{align*}
In this manner, observing the process
$n\to n+1$, we shall show that
$$
F_{n} \cap [\Mex F_{n},\Mex F_{n} + \Mex D_{n})$$
receives $v_2(n+1)$,
$$
F_{n} \cap [\Mex F_{n} + \Mex D_{n},\infty)$$
receives $v_3(n+1)$, 
and
$$
D_{n}\cap [\Mex D_{n}, \infty)
$$
receives $v_3(n+1)-v_1(n+1)$ 
 and the minimum of them may stay the same or increase by the additional 
elements. All three cases which will appear 
in the proofs of Theorem \ref{GAP} and \ref{Fin} are found 
in the above example.

For $K=\{k_i\ |\ i=1,\dots,\ell\}\subset \N_0$ with $k_{i+1}>k_i$, we define
$$
\Gap(K):=
\begin{cases} \min\{ k_{i+1}-k_i \ |\ i=1,\dots,\ell\}& \text{Card}(K)\ge 2\\ 
\infty & \text{Card}(K)\le 1.
\end{cases}                            
$$
The dynamics on $
F_n \cap [\Mex F_n + \Mex D_n,\infty)$ 
and
$
D_n\cap [\Mex D_n, \infty)
$
are pretty easy.

\begin{thm}
\label{GAP}
We have

\begin{enumerate}[a\textrm{[n]:}]
\item $v_i(n)-v_i(n-1)\ge i$ for $i=1,2,3$
%\item $v_1(n)-2 v_2(n)+v_3(n)\in \{0,1,2\}$
\item $v_3(n)-v_1(n)\ge v_3(n-1)-v_1(n-1)+2$
\item $\Gap(D_{n-1}\cap [\Mex(D_{n-1}), \infty])\ge 2$
\item $\Gap(F_{n-1}\cap [\Mex(F_{n-1})+\Mex(D_{n-1}), \infty])\ge 3$
\item $\max F_n=v_3(n)=v_2(n)+\Mex(D_{n-1})$
\item $\max D_n=v_3(n)-v_1(n)$
\end{enumerate}
\end{thm}

\begin{proof}
$v_1(n)-v_1(n-1)\ge 1$ is clear by the definition of $\Mex$.
We prove all other statements 
%except $b[n]$ 
by induction on $n$. 
Assume that statements are valid for $n$ and less. Set $h_n=\Mex(D_n)$.
%We have $v_3(n)=v_2(n)+h_{n-1}$ by induction assumption. 

{\bf Case 1.} 
If $v_1(n)+h_{n-1}\not \in F_{n-1}$, then $v_2(n)=v_1(n)+h_{n-1}$.
Thus $h_{n-1}\in D_n$ and this implies $h_{n}\ge h_{n-1}+1$.
Then we have
\begin{equation*}
%\label{1_2}
v_2(n+1)\ge v_1(n+1)+h_n\ge v_1(n)+1+h_{n-1}+1=v_2(n)+2
\end{equation*}
and
\begin{equation}
\label{1_3}
v_3(n+1)\ge v_2(n+1)+h_n\ge v_2(n)+2+h_{n-1}+1=v_3(n)+3
\end{equation}
by $e[n]$, which proves $a[n+1]$. $e[n]$ implies
\begin{equation}
\label{31_1}
v_3(n+1)-v_1(n+1)-(v_3(n)-v_1(n))\ge  2(h_n-h_{n-1})\ge 2,
\end{equation}
which shows $b[n+1]$ and $f[n+1]$. 
$e[n]$ and (\ref{1_3}) implies $v_3(n+1)= v_2(n+1)+h_n$ and $e[n+1]$.
$b[n+1]$ and $h_n>h_{n-1}$ gives $c[n+1]$ and $f[n+1]$.
$$
\Mex (F_n)+\Mex (D_n)=v_1(n+1)+h_n\ge v_1(n)+h_{n-1}+2=\Mex (F_{n-1})+\Mex (D_{n-1})+2
$$
and $a[n+1]$ gives $d[n+1]$.

{\bf Case 2.}
If $v_1(n)+h_{n-1}\in F_{n-1}$ and $h_{n-1}+1\not \in D_{n-1}$, then
$d[n]$ gives $v_1(n)+h_{n-1}+1\not \in F_{n-1}$ and therefore $v_2(n)=v_1(n)+h_{n-1}+1$.
Since $v_3(n)=v_2(n)+h_{n-1}$ by $e[n]$, we have $h_{n-1},h_{n-1}+1\in D_{n}$
and $h_{n}\ge h_{n-1}+2$.
Then we have
\begin{equation*}
%\label{2_2}
v_2(n+1)\ge v_1(n+1)+h_n\ge v_1(n)+1+h_{n-1}+2=v_2(n)+2
\end{equation*}
and
\begin{equation}
\label{2_3}
v_3(n+1)\ge v_2(n+1)+h_n\ge v_2(n)+2+h_{n-1}+2=v_3(n)+4
\end{equation}
by $e[n]$, which proves $a[n+1]$.  $e[n]$ implies
\begin{equation}
\label{31_2}
v_3(n+1)-v_1(n+1)-(v_3(n)-v_1(n))\ge  2h_n-(2h_{n-1}+1)\ge 3,
\end{equation}
which shows $b[n+1]$ and $f[n+1]$.
$e[n]$ and (\ref{2_3}) implies $v_3(n+1)= v_2(n+1)+h_n$ and $e[n+1]$.
$b[n+1]$ and $h_n>h_{n-1}+1$ gives $c[n+1]$.
$$
\Mex (F_n)+\Mex (D_n)=v_1(n+1)+h_n\ge v_1(n)+h_{n-1}+3=\Mex (F_{n-1})+\Mex (D_{n-1})+3
$$
and $a[n+1]$ gives $d[n+1]$.

{\bf Case 3.}
If $v_1(n)+h_{n-1}\in F_{n-1}$ and $h_{n-1}+1\in D_{n-1}$, then
$d[n]$ gives $v_1(n)+h_{n-1}+2\not \in F_{n-1}$ and therefore $v_2(n)=v_1(n)+h_{n-1}+2$.
Since $v_3(n)=v_2(n)+h_{n-1}$ by $e[n]$, we have $h_{n-1},h_{n-1}+1,h_{n-1}+2\in D_{n}$
and $h_{n}\ge h_{n-1}+3$.
Then we have
\begin{equation*}
%\label{3_2}
v_2(n+1)\ge v_1(n+1)+h_n\ge v_1(n)+1+h_{n-1}+3=v_2(n)+2
\end{equation*}
and
\begin{equation}
\label{3_3}
v_3(n+1)\ge v_2(n+1)+h_n\ge v_2(n)+2+h_{n-1}+3=v_3(n)+5
\end{equation}
by $e[n]$,
which proves $a[n+1]$.  $e[n]$ implies
\begin{equation}
\label{31_3}
v_3(n+1)-v_1(n+1)-(v_3(n)-v_1(n))\ge  2h_n-(2h_{n-1}+2)\ge 4,
\end{equation}
which shows $b[n+1]$ and $f[n+1]$.
$e[n]$ and (\ref{3_3}) implies $v_3(n+1)= v_2(n+1)+h_n$ and $e[n+1]$.
$b[n+1]$ and $h_n>h_{n-1}+2$ gives $c[n+1]$.
$$
\Mex (F_n)+\Mex (D_n)=v_1(n+1)+h_n\ge v_1(n)+h_{n-1}+4=\Mex (F_{n-1})+\Mex (D_{n-1})+4
$$
and $a[n+1]$ gives $d[n+1]$.
\end{proof}

The statement $b[n]$ is generalized to

\begin{cor}
\label{2Diff}
For $i>j$, we have
$$
v_i(n+1)-v_j(n+1)-(v_i(n)-v_j(n)) \ge i-j
$$
\end{cor}

\begin{proof}
$b[n]$ is for $i=3, j=1$.  $e[n]$ implies
$$
v_3(n+1)-v_2(n+1)-(v_3(n)-v_2(n)) = h_{n}-h_{n-1}\ge 1.
$$
For $i=2,j=1$, we follow the proof of Theorem \ref{GAP}. In Case 1, 
we have
$$
v_2(n+1)-v_1(n+1)-(v_2(n)-v_1(n)) \ge h_{n}-h_{n-1}\ge 1.
$$
Further in Case 2, 
$$
v_2(n+1)-v_1(n+1)-(v_2(n)-v_1(n)) \ge h_{n}-h_{n-1}-1\ge 1,
$$
and in Case 3, 
$$
v_2(n+1)-v_1(n+1)-(v_2(n)-v_1(n)) \ge h_{n}-h_{n-1}-2\ge 1.
$$
\end{proof}

\begin{cor}
\label{Balance}
$v_1(n)-2 v_2(n)+v_3(n)\in \{0,-1,-2\}$.
\end{cor}

\begin{proof}
In the proof of Theorem \ref{GAP}, we showed 
$v_2(n)-v_1(n)-h_{n-1}\in \{0,1,2\}$. By $e[n]$, we also have
 $v_3(n)-v_2(n)-h_{n-1}=0$. 
\end{proof}

From the proof of Theorem \ref{GAP}, 
we obtain an understanding of the process to determine
the sequence $(v_1(n),v_2(n),v_3(n))$. Define 
$$
\kappa(n):= \begin{cases} 0 & \Mex(F_n)+\Mex(D_{n})\not\in F_{n}\\
1 & \Mex(F_n)+\Mex(D_{n})\in F_{n}\ \land\ 1+\Mex(D_n)\not\in D_n\\
2 & \Mex(F_n)+\Mex(D_{n})\in F_{n}\ \land\ 1+\Mex(D_n)\in D_n.
\end{cases}
$$
From the three cases of the proof, we see
$$
F_{n+1}=F_n\cup \{\Mex(F_n),\Mex(F_n)+\Mex(D_{n})+\kappa(n), 
\Mex(F_n)+2\Mex(D_{n})+\kappa(n)\}$$ 
and
$$
D_{n+1}=D_n\cup \{\Mex(D_{n}),\Mex(D_{n})+\kappa(n),2\Mex(D_{n})+\kappa(n)\}.
$$
Note that new elements inductively
fill the smallest lacuna of 
$F_n$ and $D_n$. Therefore the beginning of $F_n$ and $D_n$ are 
consecutive integers and length of consecutive integers
steadily increases. 
We only have to remember $(\Mex(F_n),\Mex(D_n),W_n,Z_n)$
where
$$
W_n=F_n \cap [\Mex(F_n),\infty), \ Z_n=D_n \cap [\Mex(D_n),\infty).
$$
Summarizing these, we obtain a simplified algorithm:

\begin{itemize}
\item $M=1, D=1, W=\{\}, Z=\{\},S=\{(0,0,0)\}$.
\item[\dag] $\kappa:= \begin{cases} 0 & M+D \not\in W\\
1 & M+D \in W\ \land\ D+1\not\in Z\\
2 & M+D \in W\ \land\ D+1\in Z
\end{cases}$
\item $S\leftarrow S\cup \{(M,M+D+\kappa,M+2D+\kappa)\}$
\item $W\leftarrow W \cup \{M,M+D+\kappa,M+2D+\kappa\}$
\item $Z\leftarrow Z\cup \{D, D+\kappa,2D+\kappa\}$
\item $M\leftarrow \min((\N\cap [M,\infty)) \setminus W), D\leftarrow 
\min((\N \cap [D,\infty)) \setminus Z)$
\item $W\leftarrow W\cap [M, \infty),\ Z\leftarrow Z\cap [D,\infty)$
\item Go to $\dag$
\end{itemize}
which output $(\Mex(F_n),\Mex(D_n),W_n,Z_n)\ (n=1,2,\dots)$
and the set of P-positions 
$\{(v_1(j),v_2(j),v_3(j))\ |\ j=0,1,\dots, n \}$.
In this way, the required memory of our algorithm is made small. However 
the sequence may not be easy, 
since the set $W_n$, $Z_n$ can have unbounded cardinality.
It is illuminating to introduce
$$
X_n=W_n\cap [\Mex(F_n),\Mex(F_n)+\Mex(D_n)),\quad 
Y_n=W_n\cap [\Mex(F_n)+\Mex(D_n),\infty).
$$
By the inductive process described above, $X_n$ receives $v_2(n+1)$ and
$Y_n$ receives $v_3(n+1)$, and the minimum of $X_n$ increases by the addition of
$v_1(n+1)$. The same is true for $Z_n$. It receives the maximum element
$v_3(n+1)-v_1(n+1)$ and the minimum increases (or equal) by the new elements
$\{\Mex(D_n),\Mex(D_n)+\kappa(n)\}$. By $d[n]$, the gaps in 
$Y_n$ is greater or equal to $3$, and by $c[n]$ the gaps in $Z_n$ is greater or equal to $2$. 

\section{Finite Difference Property}
\label{FD}

For $K=\{k_i\ |\ i=1,\dots,\ell\}\subset \N_0$ with $k_{i+1}>k_i$, let 
$$\Bl(K)=\max\{j\ |\ k_{i+1}-k_{i}=k_{i+2}-k_{i+1}=\dots=k_{i+j-1}-k_{i+j-2}=
1\},$$
i.e., the maximum length of arithmetic progression of difference $1$ in $K$.
We wish to study the dynamics of $X_n$
as $n\to n+1$.

\begin{thm}
\label{Fin}
We have
\begin{align}
\label{Run}
\Bl(F_n \cap [\Mex(F_n),\Mex(F_n)+&\Mex(D_n)-1])\le 3\\
\label{z1}
v_1(n+1)-v_1(n)&\in \{1,2,3,4\},\\
\label{z2}
v_2(n+1)-v_2(n)&\in \{2,3,4,5,6,7,8\},\\
\label{z3}
v_3(n+1)-v_3(n)&\in \{3,4,5,6,7,8,9,10,11,12\}\\
\label{d}
\Mex(D_n)-\Mex(D_{n-1})&\in \{1,2,3,4\}.
\end{align}

\end{thm}

\begin{proof}
Note that $$F_{n} \cap [1,\Mex(F_{n})-1]=\N_0\cap  [1,\Mex(F_n)-1]$$ 
and 
$$\Gap(F_{n}\cap [\Mex(F_{n})+\Mex(D_{n}), \infty])\ge 3$$
by $d[n+1]$. We are interested in the middle range
$X_n$. As $n\to n+1$, we move to
$X_{n+1}$ and the both ends are shifted to the right. 
We add three points $v_1(n+1),v_2(n+1),v_3(n+1)$ to $F_n$. 
To $X_n$, we simply 
add a new point 
$v_2(n+1)=v_1(n+1)+\Mex(D_n)+i$ with $i\in \{0,1,2\}$.
We are interested in this iteration as $n\to \infty$.

Recall that addition of $v_2(n+1)$ occurs at the lacuna of $Y_n$
where $\Gap(Y_n)\ge 3$. 
Since $v_2(n+1)-v_2(n)\ge 2$ by $a[n]$, 
it is impossible to produce $4$ consecutive points in $X_n$
for any $n$, 
in the iteration of insertion of elements of $v_2$ values. 
We have shown (\ref{Run})
which immediately implies (\ref{z1}).

We may assume that 
$h_{n}=\Mex(D_n)\ge 3$. Recalling 
the proof of Theorem \ref{GAP}, we have

{\bf Case 1 $v_2(n)=v_1(n)+h_{n-1}$.} We have $h_{n-1}\in D_n$ and
$D_n = D_{n-1} \cup \{ h_{n-1},2h_{n-1} \}$. 
If $h_{n-1}+1\not \in D_{n-1}$, then since $h_{n-1}+1\not \in D_n$
we have $h_n=\Mex(D_n)=h_{n-1}+1$.
If  $h_{n-1}+1\in D_{n-1}$ then by $c[n]$, we see $h_{n-1}+2 \not\in D_{n-1}$
and therefore $h_n=\Mex(D_n)=h_{n-1}+2$.

{\bf Case 2 $v_2(n)=v_1(n)+h_{n-1}+1$.}  We have $h_{n-1}, h_{n-1}+1\in D_n$ and
$D_n = D_{n-1} \cup \{ h_{n-1},h_{n-1}+1,2h_{n-1}+1 \}$. 
If $h_{n-1}+2\not \in D_{n-1}$, then since $h_{n-1}+2\not \in D_n$
we have $h_n=\Mex(D_n)=h_{n-1}+2$.
If  $h_{n-1}+2 \in D_{n-1}$ then by $c[n]$, we see $h_{n-1}+3 \not\in D_{n-1}$
and therefore $h_n=\Mex(D_n)=h_{n-1}+3$.

{\bf Case 3 $v_2(n)=v_1(n)+h_{n-1}+2$.}  We have $h_{n-1}, h_{n-1}+1,
h_{n-1}+2\in D_n$ and
$D_n = D_{n-1} \cup \{ h_{n-1},h_{n-1}+2,2h_{n-1}+2 \}$. 
If $h_{n-1}+3\not \in D_{n-1}$, then since $h_{n-1}+3\not \in D_n$
we have $h_n=\Mex(D_n)=h_{n-1}+3$.
If  $h_{n-1}+3 \in D_{n-1}$ then by $c[n]$, we see $h_{n-1}+4 \not\in D_{n-1}$
and therefore $h_n=\Mex(D_n)=h_{n-1}+4$.

In all three cases, 
we have seen $h_n-h_{n-1}\in \{1,2,3,4\}$ which is (\ref{d}). 
We have 
\begin{align*}
v_2(n)&=v_1(n)+h_{n-1} + k\quad (k\in \{0,1,2\})\\
v_2(n+1)&=v_1(n+1)+h_n + j\quad (j\in \{0,1,2\})
\end{align*}

which implies
\begin{equation}
\label{21diff}
v_2(n+1)-v_2(n)=v_1(n+1)-v_1(n)+h_n-h_{n-1}+j-k \le 10.
\end{equation}

We can improve this bound to $8$ since in 
the above three cases, if $h_n-h_{n-1}=4$ then $k=2$
and if $h_n-h_{n-1}=3$ then $k\ge 1$. Therefore
we have $h_n-h_{n-1}+j-k\le 4$ and (\ref{z2}).
Finally

\begin{align*}
v_3(n)&=v_2(n)+h_{n-1}\\
v_3(n+1)&=v_2(n+1)+h_n 
\end{align*}

gives the bound 
\begin{equation}
\label{32diff}
v_3(n+1)-v_3(n)=v_2(n+1)-v_2(n)+h_n-h_{n-1} \le 12.
\end{equation}
and we obtain (\ref{z3}).
\end{proof}

\begin{figure}[h]
\includegraphics[clip]{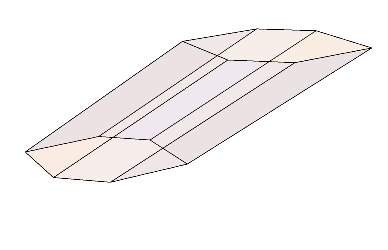}
\caption{Differences form an oblique hexagonal prism\label{Prism}}
\end{figure}

There are 56 differences of adjacent P-positions and they 
form the shape of an oblique hexagonal prism 
in Figure \ref{Prism}. This set 
gives a tiling of $\Z^3$ by lattice translations. 

\begin{cor}
\label{FinDiff}
\begin{align*}
&\{
(v_1(n+1)-v_1(n),
v_2(n+1)-v_2(n),
v_3(n+1)-v_3(n))\ |\ n\in \N_0 \}\\
&=\{ (a,b,c)\in \N^3\ |\ a\le 4,\ 
0<b-a\le 4,\ 0<c-b\le 4 ,\ |a-2b+c|\le 2\}\\
&=\{ (i,i+j,i+j+k)\ |\ i,j,k\in \{1,2,3,4\},\ |j-k|\le 2\}.
\end{align*}
\end{cor}

\begin{proof}
Last two expressions are equivalent. Theorem \ref{Algo}
implies $a,b,c\in \N$. (\ref{z1}) shows $a\le 4$.
Corollary \ref{Balance} gives 
$|a-2b+c|\le 2$. 
Corollary \ref{2Diff} and 
(\ref{21diff}) with $h_n-h_{n-1}+j-k\le 4$ implies $0<b-a\le 4$, and
 (\ref{32diff}) and (\ref{d}) 
implies $0<c-b\le 4$. Thus the first set is
a subset of the second. We show the converse inclusion by direct computation.
We write
$$
\{ (i,i+j,i+j+k)\ |\ i,j,k\in \{1,2,3,4\},\ |j-k|\le 2\}=\{ k_{\ell}\ |\ \ell=1,2,\dots, 56\}
$$
with $k_{\ell}\ll k_{\ell+1}$ by lexicographic ordering. Let $n_{\ell}$ be the minimum $n$
that
$$
k_{\ell}=(v_1(n)-v_1(n-1),
v_2(n)-v_2(n-1),
v_3(n)-v_3(n-1)).
$$
Then we have
\begin{align*}
&( n_1,n_2, \dots, n_{56})\\
&=(1, 58, 49, 11, 10, 21, 15, 31, 43, 18, 712, 48, 327, 1913, 27, 30, 19, 55, 5, 44, 90, \\
& 28, 286, 99, 1419, 89, 14, 325, 8, 155, 200, 20, 2, 12, 131, 231, 59, 384, 
686, 338, 593, \\
& 1269, 65, 445, 706, 869, 6, 150, 51, 17, 573, 4656, 939, 2724, 14449, 26185).
\end{align*}
For example $v(26184)=(49223, 93770, 138315)$ and
$v(26185)=(49227, 93778, 138327)$ yields the difference $n_{56}=(4,8,12)$
for the first time. Subsequently $(4,8,12)$ 
occurs at $n=43153, 46142, 46704, 48335, 93299, \dots$. 
\end{proof}

\section{Linear Structure}
\label{LS}

As in Figure \ref{P}, the P-positions approximate half-lines emanating from the origin:

\begin{figure}[h]
\subfigure{
\includegraphics[width=0.4\columnwidth]{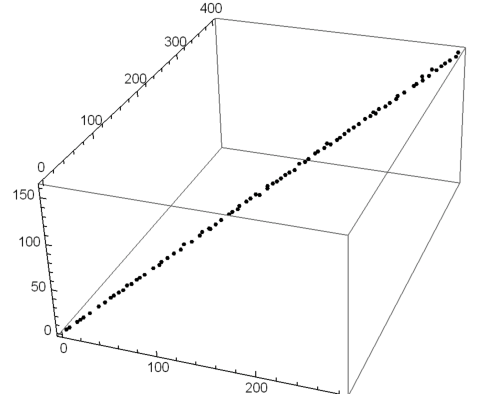}}
\subfigure{
\includegraphics[width=0.4\columnwidth]{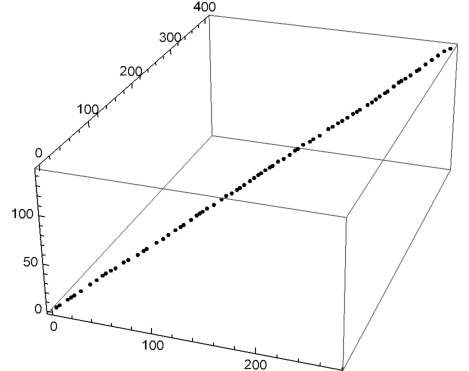}}
\caption{P-positions of positive Nim and $L^1$-Nim\label{P}}
\end{figure}

\begin{conj}
The limit
$$
(\alpha,\beta,\gamma):=
\lim_{n\to \infty}
\left(\frac{v_1(n)}n,\frac {v_2(n)}n,\frac {v_3(n)}n\right)
$$
exists.
\end{conj}
This line is possibly approximated by
$$
\frac X{\alpha} =\frac Y{\beta}=\frac Z{\gamma}.
$$
Numerical experiments suggest
$$
\alpha\approx 1.97,\quad
\beta\approx 3.49,\quad
\gamma\approx 4.88
$$
for positive Nim and 
$$
\alpha\approx 1.88,\quad
\beta\approx 3.58,\quad
\gamma\approx 5.28
$$
for $L^1$-Nim.

If the conjecture is valid for $L^1$-Nim, then from Corollary \ref{Par} we see
\begin{equation}
\label{Freq}
\frac 1{\alpha} +\frac 1{\beta} + \frac 1{\gamma}=1.
\end{equation}

Corollary \ref{Balance} gives another constraint:
\begin{equation}
\label{Arith3}
\alpha-2\beta+\gamma=0.
\end{equation}

We wish to have an additional relation among $\alpha,\beta,\gamma$ to identify
these numbers. We tried to find the
substitutive structure of the sequences $v_i(n+1)-v_i(n)$ by taking 
successive derived sequences but in vain, 
see \cite{Durand:98,Holton-Zamboni:99}.
\bigskip

%Define $\alpha_n=v_1(n)/n$,  $\beta_n=v_2(n)/n$, $\gamma_n=v_3(n)/n$. By Corollary \ref{Par}, we have $\alpha_n+\beta_n+\gamma_n=1$ and Corollary \ref{Balance}
%gives
%$$
%\alpha_n-2\beta_n+\gamma_n=O(1/n).
%$$

Let us assume that 
$\lim_{n\to \infty}(v_1(n)/n,v_2(n)/n,v_3(n)/n) 
=(\alpha,\beta,\gamma)$ exists
and put $\delta=\gamma-\beta=\beta-\alpha$.
We can rewrite them as
$$
\alpha=\beta\left(1-\sqrt{\frac {\beta-3}{\beta-1}}\right),\ \gamma=\beta\left(1+\sqrt{\frac {\beta-3}{\beta-1}}\right),\ \delta=\beta \sqrt{\frac {\beta-3}{\beta-1}}.
$$
By (\ref{d}) or Corollary \ref{FinDiff}, we have
$$
1\le \delta \le 4.
$$
Parameterizing these by a new parameter $\xi=\sqrt{\frac {\beta-3}{\beta-1}}$,
we can solve this inequality on $\beta$. We have
$$
3.21432\le \beta\le 5.41061
$$
and consequently 
$$
1.41061\le \alpha \le 2.21432, \ 
4.21432\le \gamma \le 9.41061.
$$
These $6$ numerical constants can be replaced by cubic numbers. 
Hereafter we give a heuristic discussion to determine $\alpha,\beta,\gamma$. 

Recalling the proof of Theorem \ref{GAP}, the set
$
Y_m
$
is the union of elements of $v_3(n)$, and it 
natural density $1/\gamma$ as $m\to \infty$. Similarly the set
$
Z_{m}$
is a union of $v_3(n)-v_1(n)$ and 
its natural density is $1/(2\delta)=1/(\gamma-\alpha)$.
Since we can not find any further 
structure, we assume that
an integer in $$[\Mex(F_m)+\Mex(D_m),\max(F_m)+\max(D_m))$$
falls into $F_{m}$ with probability $1/\gamma$. 
Avoiding the probabilistic terminology, we assume that
\begin{equation}
\label{Assumption}
\lim_{N\to \infty} 
\frac 1N \mathrm{Card}\left\{ n\le N \ \left|\ v_1(n+1)+\Mex(D_n)\in F_n\right.\right\}=
\frac 1{\gamma},
\end{equation}
which is equivalent to
\begin{equation}
\label{Assumption2}
\lim_{N\to \infty} 
\frac 1N \mathrm{Card}\left\{ n\le N \ \left|\ \kappa(n)=0\right.\right\}=1-
\frac 1{\gamma}.
\end{equation}
Recall that $\kappa(n)=0$ is also equivalent to $v_1(n)-2v_2(n)-v_3(n)=0$.

Let us fix $m$ and study the dynamics on
$
Z_m=D_{m}\cap [\Mex(D_{m}),\infty]
$
for $n=m,m+1,\dots$. 
In Case 1, $h_n=h_{n-1}+1, 2h_n
\in D_n$ are the additional elements when $n\to n+1$.
In Case 2, $h_n, h_n+1, 2h_n+1$ are added, and in Case 3, $h_n, h_n+2, 2h_n+2$ 
are added when $n\to n+1$.
The addition is done in a way
to fill the lacuna of $W_m$ from the left
and therefore the left side
of the set $W_m$ is gradually filled by 
new elements to form consecutive integers. 
The elements $2h_n+j$ for $j=0,1,2$ is 
located far from the left and these are not related to this discussion. 
This dynamics increases $\Mex(D_n)$ by $1,2,3,4$. 
In particular, we observe that
if $v_1(n+1)+\Mex(D_n)\not\in F_n$ then one lacuna is filled, and if
 $v_1(n+1)+\Mex(D_n)\in F_n$ then two lacuna are filled. 
By our assumption, 
$v_1(n+1)+\Mex(D_n)\not\in F_n$ happens with natural density $1-1/\gamma$
and 
$v_1(n+1)+\Mex(D_n)\in F_n$ 
happens with natural density $1/\gamma$, 
the average 
of the number of lacuna 
to be filled by the increase of $n$ is
$$
1-1/\gamma+2/\gamma=1+1/\gamma.
$$
Since $Z_m$ has density $1/(2\delta)$, the set of lacuna 
has natural density $1-1/(2\delta)$. 
After $n$-times increments, $n(1+1/\gamma)$
lacuna are filled. Considering the density of the lacuna, we see
$$
n(1+1/\gamma)=(1-1/(2\delta)) (\Mex(D_{m+n-1})-\Mex(D_m))+o(n)
$$
for any fixed $m\in \N$. In the light of $e[n]$, 
dividing by $n$ and taking $n\to \infty$, 
we see 
$$
(1+1/\gamma)=(1-1/(2\delta)) \delta
$$
or
\begin{equation}
\label{incD}
3+2/\gamma=2\delta.
\end{equation}
In this manner, we can obtain the expected values of $\alpha, \beta, \gamma$.

\begin{thm}
\label{Line}
Assume that 
$\lim_{n\to \infty}(v_1(n)/n,v_2(n)/n,v_3(n)/n) 
=(\alpha,\beta,\gamma)$ exists and the asymptotic formula 
(\ref{Assumption}) holds. 
Then
$\alpha,\beta,\gamma$ are algebraic numbers of degree $5$,
approximately
$$
(\alpha,\beta,\gamma)\approx
(1.88542779715, 3.57535140026, 5.26527500337).
$$
\end{thm}

\begin{proof}
From three relations (\ref{Freq}),(\ref{Arith3}) and (\ref{incD}), with the 
help of the parameter $\xi$, 
we obtain totally real quintic minimum polynomial
\begin{equation}
\label{MinP}
8 x^5-44 x^4+54 x^3+39 x^2-60 x-22
\end{equation}
of $\delta$ whose Galois group is the symmetric group\footnote{
The minimum polynomial of $2\delta $ is factored into
$\left(x^2+8 x+11\right) \left(x^3+15 x^2+15 x+9\right) \pmod{17}$
where $17$ is not a divisor of its discriminant. This shows 
the existence of a transposition in the Galois group.} $S_5$.
Thus $\delta$ is not solvable by radicals over $\Q$.
Under the condition $\beta>3$, the 
solution is unique and we have
$$
(\alpha,\beta,\gamma,\delta)\approx
(1.88542779715, 3.57535140026, 5.26527500337, 1.68992360311)
$$
Here
$$
\alpha
=2 \delta^4-8 \delta^3+\frac{3 \delta^2}{2}+10 \delta+3
$$
and $\beta=\alpha+\delta,\ \gamma=\alpha+2\delta$.
\end{proof}

\begin{figure}[h]
\includegraphics{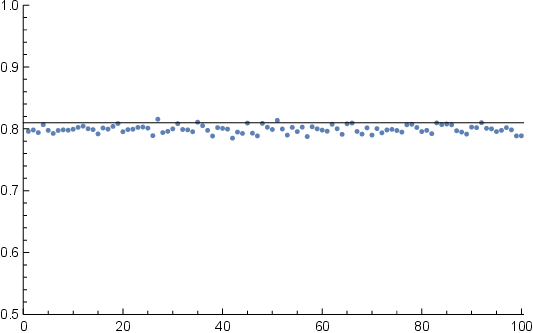}
\caption{The frequency of $n\in [5000j, 5000(j+1)-1]$
which satisfy $v_1(n)-2v_2(n)+v_3(n)=0$ for $j=0,1,\dots, 
100.$\label{FreqKappa}}
\end{figure}

Finally we discuss the adequacy of our assumption (\ref{Assumption}).
Recall that the sequence $\Mex(F_n)+\Mex(D_n)$ is expected to 
increase by $\beta$ and $Y_n$ consists of $v_3(n)$
which is expected to increase by $\gamma$ as $n\to n+1$.
Since $1, \beta$ and $\gamma$ are linearly independent over $\Q$, 
two sequences $n\beta\ (n\in \Z)$ and $m\gamma\ (m\in \Z)$ 
have no linear relation, and our probabilistic assumption does not collapse.
%Figure \ref{FreqKappa} gives the result of a
A statistical experiment for (\ref{Assumption2})
 is done: the frequency of $n\in [1,500,000]$ with $\kappa(n)=0$
is $0.7993\dots$. Together with Figure \ref{FreqKappa}, 
the assumption (\ref{Assumption2}) might be acceptable
with $1-1/\gamma\approx 0.810076$.
Further, we have a good numerical fit:
$$
1-\left(\frac 1{m}\sum_{n=1}^{m} \frac{v_3(n)}n \right)^{-1}
\approx 0.810706 \qquad \text{ with }  m=500,000.
$$
It is of interest to make experiments in the larger range.

\bigskip

{\bf Acknowledgments}
\medskip

I am indebted to Bo Tan who suggested me this $L^1$ rule
when I stayed in Wuhan China on March 2015.
I also would like to thank David Klein for the
discussion and 
sending me the preprint \cite{KleinFraenkel} during 2022.
This work had been my long project without official output
and finally came into this
shape after I posed this game at
the `Open problem workshop" in Hayama on September 2024.
I wish to show my gratitude to the participants
 who shared this problem there.

%\bibliographystyle{amsplain}
%\bibliography{../../reflist}

\begin{thebibliography}{1}

\bibitem{KleinFraenkel}
D.~Klein and A.S.~Fraenkel, \emph{K-pile {W}ythoff games, {L}exicographic
  {W}ythoff}, to appear in Games Of No Chance.

\bibitem{WythoffV}
E.~Duch\^ene, A.S. Fraenkel, V.~Gurvich, N.B. Ho, C.~Kimberling, and
  U.~Larsson, \emph{Wythoff visions}, Games of no chance 5, Math. Sci. Res.
  Inst. Publ., vol.~70, Cambridge Univ. Press, Cambridge, 2019, pp.~35--87.

\bibitem{Durand:98}
F.~Durand, \emph{A characterization of substitutive sequences using return
  words}, Discrete Math. \textbf{179} (1998), no.~1-3, 89--101.

\bibitem{Holton-Zamboni:99}
C.~Holton and L.Q. Zamboni, \emph{Descendants of primitive substitutions},
  Theory Comput. Syst. \textbf{32} (1999), no.~2, 133--157.

\bibitem{Rigo_Game}
M.~Rigo, \emph{From combinatorial games to shape-symmetric morphisms},
  Substitution and tiling dynamics: introduction to self-inducing structures,
  Lecture Notes in Math., vol. 2273, Springer, Cham, [2020] \copyright 2020,
  pp.~227--291.

\bibitem{Wythoff}
W.A. Wythoff, \emph{A modification of the game of {N}im}, Nieuw. Arch. Wiskd.
  \textbf{7} (1907), 199--202.

\end{thebibliography}

\end{document}